\documentclass[table]{amsart} 

\usepackage{euscript}           % Nice script font.                       
\usepackage{pstricks}
\usepackage{pst-plot}
\usepackage{float,lscape}
\input{multido.tex}
\usepackage{xcolor}

\usepackage{pdfpages}

\usepackage{amsrefs,amsmath,amsthm,amssymb}            
\usepackage{graphicx,enumerate,calc,lscape}
\usepackage[matrix,arrow,curve,frame]{xy}   

\usepackage{longtable,booktabs}

\usepackage{amsaddr}

\newtheorem{thm}[subsection]{Theorem}

\theoremstyle{definition}  % Bold headings and Roman body text.

\newcommand{\map}{\rightarrow}
\newcommand{\tmf}{\textit{tmf}}
\newcommand{\mmf}{\textit{mmf}}

\newcommand{\field}[1]  {\mathbb #1} % Use blackboard bold for these sets

\newcommand{\F}         {\field F}

\newcommand{\C}         {\field C}
\newcommand{\M}         {\field M}

\newcommand{\D}{\Delta}

\DeclareMathOperator{\Sq}{Sq}

\DeclareMathOperator{\Ext}{Ext}

\numberwithin{equation}{subsection}

\newrgbcolor{grey}{0.5 0.5 0.5}

\begin{document}

\title{The homotopy of $\C$-motivic modular forms}

\author{Daniel C.\ Isaksen}
\address{Department of Mathematics \\ Wayne State University\\
Detroit, MI 48202}
\email{isaksen@wayne.edu}

%%%%%%%%%%%%%%%%%%%%%%%%%%%%%%%%%%%%%%%%%%%%%%%%%%%%%%%%%%%%%%%%%%

\begin{abstract}
A $\C$-motivic modular forms spectrum $\mmf$ has recently been
constructed.  
This article presents detailed computational information on the
Adams spectral sequence for $\mmf$.
This information is essential for computing
with the $\C$-motivic and classical Adams spectral sequences
that compute the $\C$-motivic and classical stable homotopy
groups of spheres.
\end{abstract}

\maketitle

%%%%%%%%%%%%%%%%%%%%%%%%%%%%%%%%%%%%%%%%%%%%%%%%%%%%%%%%%%%%%%%%%%

\section{Introduction}

\subsection{Topological modular forms}
The topological modular forms spectrum $\tmf$ is an
important tool for studying classical stable 
homotopy groups at the prime $2$.
Its unit map $S^0 \map \tmf$ induces a map of Adams spectral sequences,
and much information about the behavior of the Adams spectral sequence
for $S^0$ can be deduced from the behavior of the Adams spectral
sequence for $\tmf$.

From this perspective, the essential properties of $\tmf$
are that it is a ring spectrum whose $\F_2$-cohomology is
$A//A(2)$.  Here $A$ is the usual Steenrod algebra,
and $A(2)$ is the subalgebra generated by $\Sq^1$, $\Sq^2$, and $\Sq^4$.
A standard change-of-rings theorem then implies that the 
$E_2$-page of the Adams spectral sequence for $\tmf$ is the
cohomology of $A(2)$, i.e., $\Ext_{A(2)} (\F_2, \F_2)$.

The cohomology of $A(2)$ was first computed by May in his thesis \cite{May64}.
The computation is also described in 
\cite{BR} 
\cite{tmf14}*{Chapter 13}  
\cite{Isaksen09}
\cite{SI67}.
It turns out that the Adams differentials are manageable,
and a complete description of the $E_\infty$-page is possible.
It is also possible to determine all of the hidden extensions
by $2$, $\eta$, and $\nu$.  As a result, we thoroughly understand
the ring of stable homotopy groups of $\tmf$.

The spectrum $\tmf$ has also been completely analyzed from the
perspective of the Adams-Novikov spectral sequence \cite{Bauer08}.

\subsection{Motivic modular forms}
A similar story plays out in 
$\C$-motivic stable homotopy theory.  
A $\C$-motivic modular forms spectrum $\mmf$ has recently
been constructed \cite{GIKR18}.
Specifically, $\mmf$ is a $\C$-motivic ring spectrum
whose cohomology is $A // A(2)$, where $A$ is the $\C$-motivic
Steenrod algebra and $A(2)$ is the subalgebra of $A$ generated
by $\Sq^1$, $\Sq^2$, and $\Sq^4$.  The goal is to compute
the $\C$-motivic stable homotopy groups of $\mmf$ via the 
$\C$-motivic Adams spectral sequence.
The $\C$-motivic stable homotopy groups of $\mmf$
were first computed (speculatively) by means of the
$\C$-motivic Adams-Novikov spectral sequence \cite{Isaksen09}.

The first step in our project is to identify the $E_2$-page
of the $\C$-motivic Adams spectral sequence for $\mmf$, i.e.,
the cohomology of $\C$-motivic $A(2)$.  This calculation is completely
described in \cite{Isaksen09}.

The next step is to determine Adams differentials.  This turns out to be
fairly straightforward.  The Betti realization map allows us to
deduce motivic differentials from known classical differentials.

The final step is to find hidden extensions in the
$E_\infty$-page of the $\C$-motivic Adams spectral sequence for 
$\mmf$.  As in the classical case, it is possible to resolve
all hidden extensions by $2$, $\eta$, and $\nu$, but this manuscript
does not 
provide exhaustive information about these extensions.
On the other hand, we do exhaustively list the hidden
$\tau$ extensions.

The Betti realization functor from
$\C$-motivic stable homotopy theory to classical stable homotopy theory
is well-understood from a calculational perspective \cite{Isaksen14c}.
Namely, Betti realization corresponds to inverting a particular
element $\tau$ in $\pi_{0,-1}$.
In particular, the Betti realization of $\mmf$ is $\tmf$.  Therefore,
there is a well-understood map from the $\C$-motivic Adams spectral
sequence of $\mmf$ to the classical Adams spectral sequence for $\tmf$.

In other words, the $\C$-motivic Adams spectral sequence
for $\mmf$ is compatible with the classical Adams spectral sequence
for $\tmf$.
Depending on perspective, this could mean two different things.
First, information about $\tmf$ can be transported to information
about $\mmf$.  This does not completely determine the
$\C$-motivic Adams spectral sequence for $\mmf$, 
because of the presence of $\tau$ torsion.
Second, an independent analysis of the $\C$-motivic Adams spectral
sequence for $\mmf$ leads to an analysis of the
classical Adams spectral sequence for $\tmf$ as a corollary.
This latter ``motivic-to-classical" perspective is at the heart
of recent progress in the computation of classical stable 
homotopy groups \cite{Isaksen14c}.

Without further mention, we are always working in 
$\C$-motivic stable homotopy theory that has been completed at
$2$ and at $\eta$.  This means that the Adams spectral sequence
converges.

\subsection{Logical assumptions}

We assume without further reference the entire
algebraic structure of the cohomology of
$\C$-motivic $A(2)$, as described in \cite{Isaksen09}
(especially Theorem 4.13).  

We also assume complete
knowledge of the $\C$-motivic Adams-Novikov spectral
sequence for $\mmf$, as shown in the charts in \cite{Isaksen09}.
This means that the homotopy groups of $\mmf$ are already known,
and we just have to figure out which Adams differentials
and which hidden extensions are compatible with that calculation.
See also \cite{Bauer08} for a thorough analysis of the
classical Adams-Novikov spectral sequence for $\tmf$.

We also assume complete knowledge of the classical Adams spectral
sequence for $\tmf$ \cite{tmf14}*{Chapter 13} \cite{BR}.

\subsection{Organization}

Section \ref{sctn:differentials} gives a few explicit proofs for the
handful of Adams differentials that do not follow immediately
from comparison with the classical case.

Section \ref{sctn:tables} summarizes the essential parts of the
calculation in table format.  We give multiplicative generators
for each page of the spectral sequence, as well as the values
of the Adams differentials on these generators.

\subsection{Notation}

We adopt the notation of \cite{Isaksen09} and \cite{Isaksen14c},
with the exception that the elements
$\alpha$ and $\nu$ from \cite{Isaksen09} are written here
as $a$ and $n$.

Note that the elements $a$, $n$, $\Delta h_1$, and $\Delta c$
are called $\alpha$, $\beta$, $\gamma$, and $\delta$ in \cite{BR}
and \cite{tmf14}*{Chapter 13}.

\section{Adams differentials and hidden extensions}
\label{sctn:differentials}

There are several possible approaches to computing $\C$-motivic
Adams differentials for $\mmf$:
\begin{enumerate}
\item
Compare to the known classical case of $\tmf$.  This method
takes care of most of the motivic differentials.
\item
Compare to the results of the known $\C$-motivic Adams-Novikov
spectral sequence for $\mmf$.  Since the $\C$-motivic stable
homotopy groups of $\mmf$ are already known, one can often
reverse engineer Adams differentials to obtain that answer.
\item
Use the algebraic Novikov spectral sequence \cite{Novikov67}
\cite{Miller75} that computes the $E_2$-page of the
classical Adams-Novikov spectral sequence for $\tmf$.
As discussed in \cite{GWX18}, this entirely algebraic spectral
sequence is isomorphic to the $\C$-motivic Adams spectral sequence
for $\mmf/\tau$.  Then one can pull back or push forward
differentials to the $\C$-motivic Adams spectral sequence
for $\mmf$.  While this method is extremely powerful, we will
not use it in this manuscript.
\end{enumerate}

\begin{thm}
\label{thm:E2-gen}
The multiplicative generators of the
Adams $E_2$-page for $\mmf$ are listed in
Table \ref{tab:E2-gen}, together with the 
values of the $d_2$ differential on these generators.
\end{thm}

\begin{proof}
Most of the values of the motivic $d_2$ differential
follow from the analogous classical calculations.
A few differentials require additional arguments.

The value of $d_2(e)$ follows from comparison to the
$h_1$-periodic calculation \cite{GI15}.
Then the value of $d_2(u)$ follows from the relation
$c u = h_1^2 e$.
Finally, the value of $d_2(\D u + \tau n g)$ follows 
from the relation $h_1 \cdot (\D u + \tau n g) = u \cdot \D h_1$.
\end{proof}

\begin{thm}
\label{thm:E3-gen}
The multiplicative generators of the
Adams $E_3$-page for $\mmf$ are listed in
Table \ref{tab:E3-gen}, together with the 
values of the $d_3$ differential on these generators.
\end{thm}

\begin{proof}
The structure of the $E_3$-page is completely determined by the
$d_2$ differentials in Theorem \ref{thm:E2-gen}.  Then the
$E_3$-page generators can be determined by naive algebraic
considerations.

Most of the values of the motivic $d_3$ differential
follow from the analogous classical calculations.
A few differentials require additional arguments.

There is a classical differential $d_4(a e g) = P^2 (\D c + a g)$.
This must correspond to a motivic differential
$d_4(\tau^3 a e g) = P^2 (\D c + \tau a g)$.
Therefore, $\tau^2 a e g$ cannot survive to the $E_4$-page,
and $d_3(\tau^2 a e g)$ must equal $P \D h_1^2 d$.

Essentially the same argument shows that
$d_3(\tau^2 n g^2) = \D h_1^2 d^2$, using the classical
differential $d_4(n g^2) = P (\D c + a g) d$.

From the motivic Adams-Novikov spectral sequence for $\mmf$,
the element $h_1 d g$ detects an element of
$\pi_{35,21} \mmf$ that is $\tau$-free.  Therefore,
there must be a hidden $\tau$
extension from $\tau^2 h_1 d g$ to $P a n$.
Then there is also a hidden
$\tau$ extension from $\tau^2 h_1 d g^2$ to $P a n g$.
We know from the motivic Adams-Novikov spectral sequence
for $\mmf$ that $h_1 d g^2$ detects an element of
$\pi_{55, 33}$ that is annihilated by $\tau^4$.
Therefore, $\tau P a n g$ must be hit by some differential.
The only possibility is that 
$d_3(\D^2 c)$ equals $\tau P a n g$.
\end{proof}

\begin{thm}
\label{thm:E4-gen}
The multiplicative generators of the
Adams $E_4$-page for $\mmf$ are listed in
Table \ref{tab:E4-gen}, together with the 
values of the $d_4$ differential on these generators.
\end{thm}

\begin{proof}
The structure of the $E_4$-page is completely determined by the
$d_3$ differentials in Theorem \ref{thm:E3-gen}.  Then the
$E_4$-page generators can be determined by naive algebraic
considerations.

Most of the values of the motivic $d_4$ differential
follow from the analogous classical calculations.
One differential requires further explanation.

We have a relation $P \cdot \D^6 h_1 c = \tau^4 \D^4 P e g^2$.
Therefore,
\[
P \cdot d_4(\D^6 h_1 c) = \tau^2 P g \cdot d_4 (\tau^2 \D^4 e g ) =
\tau^2 \D^4 P^2 d^2 g.
\]
It follows that $d_4(\D^6 h_1 c)$ must equal $\tau^2 \D^4 P d^2 g$.
\end{proof}

\begin{thm}
\label{thm:Einfty-gen}
The $\C$-motivic Adams spectral sequence for
$\mmf$ collapses at the $E_5$-page, i.e., the $E_5$-page
equals the $E_\infty$-page.
The multiplicative generators of the
Adams $E_\infty$-page are listed in
Table \ref{tab:Einfty-gen}.
\end{thm}

\begin{proof}
The structure of the $E_5$-page is completely determined by the
$d_4$ differentials in Theorem \ref{thm:E4-gen}.  Then the
$E_5$-page generators can be determined by naive algebraic
considerations.

For degree reasons, there are no possible $d_r$ differentials
on these generator for $r \geq 5$.
\end{proof}

Having obtained the $E_\infty$-page of the Adams spectral
sequence for $\mmf$, the final step is to determine hidden
extensions.  

\begin{thm}
\label{thm:hidden-tau}
Table \ref{tab:hidden-tau} lists some hidden extensions by $\tau$.
All hidden extensions are $g$ and $\D^8$ multiples of these.
\end{thm}

\begin{proof}
These hidden extensions are the only possibilities that are
consistent with the computations of $\pi_{*,*} \mmf$
from the motivic Adams-Novikov spectral sequence \cite{Isaksen09}.
\end{proof}

\newpage

\section{Tables}
\label{sctn:tables}

\begin{longtable}{lll}
\caption{$E_2$-page multiplicative generators
\label{tab:E2-gen} 
} \\
\toprule
$(s,f,w)$ & $x$ & $d_2(x)$ \\
\midrule \endfirsthead
\caption[]{$E_2$-page multiplicative generators} \\
\toprule
$(s,f,w)$ & $x$ & $d_2(x)$ \\
\midrule \endhead
\bottomrule \endfoot
$(0, 1, 0)$ & $h_0$ & $$ \\
$(1, 1, 1)$ & $h_1$ & $$ \\
$(3, 1, 2)$ & $h_2$ & $$ \\
$(8, 3, 5)$ & $c$ & $$ \\
$(8, 4, 4)$ & $P$ & $$ \\
$(11, 3, 7)$ & $u$ & $h_1^2 c$ \\
$(12, 3, 6)$ & $a$ & $P h_2$ \\
$(14, 4, 8)$ & $d$ & $$ \\
$(15, 3, 8)$ & $n$ & $h_0 d$ \\
$(17, 4, 10)$ & $e$ & $h_1^2 d$ \\
$(20, 4, 12)$ & $g$ & $$ \\
$(25, 5, 13)$ & $\Delta h_1$ & $$ \\
$(32, 7, 17)$ & $\Delta c + \tau a g$ & $$ \\
$(35, 7, 19)$ & $\Delta u + \tau n g$ & $h_1^2 (\Delta c + \tau a g)$ \\
$(48, 8, 24)$ & $\Delta^2$ & $\tau^2 a n g$ \\
\end{longtable}

\begin{longtable}{lll}
\caption{$E_3$-page multiplicative generators
\label{tab:E3-gen}
} \\
\toprule
$(s,f,w)$ & $x$ & $d_3(x)$ \\
\midrule \endfirsthead
\caption[]{$E_3$-page multiplicative generators} \\
\toprule
$(s,f,w)$ & $x$ & $d_3(x)$ \\
\midrule \endhead
\bottomrule \endfoot
$(0, 1, 0)$ & $h_0$ & $$ \\
$(1, 1, 1)$ & $h_1$ & $$ \\
$(3, 1, 2)$ & $h_2$ & $$ \\
$(8, 3, 5)$ & $c$ & $$ \\
$(8, 4, 4)$ & $P$ & $$ \\
$(12, 6, 6)$ & $h_0^3 a$ & $$ \\
$(14, 4, 8)$ & $d$ & $$ \\
$(17, 4, 9)$ & $\tau e$ & $P c$ \\
$(20, 4, 12)$ & $g$ & $$ \\
$(24, 6, 12)$ & $a^2$ & $\tau P h_1 d$ \\
$(25, 5, 13)$ & $\Delta h_1$ & $$ \\
$(27, 6, 14)$ & $a n$ & $$ \\
$(30, 6, 16)$ & $n^2$ & $\tau h_1 d^2$ \\
$(32, 7, 16)$ & $\tau^2 a g$ & $$ \\
$(32, 7, 17)$ & $\Delta c + \tau a g$ & $$ \\
$(36, 10, 18)$ & $h_0 a^3$ & $$ \\
$(48, 9, 24)$ & $\Delta^2 h_0$ & $$ \\
$(49, 9, 25)$ & $\Delta^2 h_1$ & $\tau^3 d^2 g$ \\
$(49, 11, 26)$ & $\tau^2 a e g$ & $P \Delta h_1^2 d$ \\
$(51, 9, 26)$ & $\Delta^2 h_2$ & $$ \\
$(55, 11, 30)$ & $\tau^2 n g^2$ & $\Delta h_1^2 d^2$ \\
$(56, 11, 29)$ & $\Delta^2 c$ & $\tau P a n g$ \\
$(60, 14, 30)$ & $\Delta^2 h_0^3 a$ & $$ \\
$(80, 15, 40)$ & $\tau \Delta^2 \Delta c$ & $$ \\
$(84, 18, 42)$ & $\Delta^2 h_0 a^3$ & $$ \\
$(96, 16, 48)$ & $\Delta^4$ & $\tau^8 n g^4$ \\
\end{longtable}

\begin{longtable}{lll}
\caption{$E_4$-page multiplicative generators
\label{tab:E4-gen}
} \\
\toprule
$(s,f,w)$ & $x$ & $d_4(x)$ \\
\midrule \endfirsthead
\caption[]{$E_4$-page multiplicative generators} \\
\toprule
$(s,f,w)$ & $x$ & $d_4(x)$ \\
\midrule \endhead
\bottomrule \endfoot
$(0, 1, 0)$ & $h_0$ & $$ \\
$(1, 1, 1)$ & $h_1$ & $$ \\
$(3, 1, 2)$ & $h_2$ & $$ \\
$(8, 3, 5)$ & $c$ & $$ \\
$(8, 4, 4)$ & $P$ & $$ \\
$(12, 6, 6)$ & $h_0^3 a$ & $$ \\
$(14, 4, 8)$ & $d$ & $$ \\
$(20, 4, 12)$ & $g$ & $$ \\
$(24, 7, 12)$ & $h_0 a^2$ & $$ \\
$(25, 5, 13)$ & $\Delta h_1$ & $$ \\
$(27, 6, 14)$ & $a n$ & $$ \\
$(31, 8, 16)$ & $\tau^2 d e$ & $P^2 d$ \\
$(32, 7, 16)$ & $\tau^2 a g$ & $$ \\
$(32, 7, 17)$ & $\Delta c + \tau a g$ & $$ \\
$(36, 10, 18)$ & $h_0 a^3$ & $$ \\
$(37, 8, 20)$ & $\tau^2 e g$ & $P d^2$ \\
$(44, 10, 22)$ & $\tau^2 a^2 g$ & $P^2 a n$ \\
$(48, 9, 24)$ & $\Delta^2 h_0$ & $\tau P \Delta h_1 d$ \\
$(49, 11, 25)$ & $\tau^3 a e g$ & $P^2 (\Delta c + \tau a g)$ \\
$(50, 10, 26)$ & $\tau^2 n^2 g$ & $P a n d$ \\
$(51, 9, 26)$ & $\Delta^2 h_2$ & $$ \\
$(55, 11, 29)$ & $\tau^3 n g^2$ & $P (\Delta c + \tau a g) d$ \\
$(56, 11, 28)$ & $\tau \Delta^2 c$ & $$ \\
$(60, 14, 30)$ & $\Delta^2 h_0^3 a$ & $$ \\
$(72, 15, 36)$ & $\Delta^2 h_0 a^2$ & $$ \\
$(80, 15, 40)$ & $\tau \Delta^2 \Delta c$ & $$ \\
$(84, 18, 42)$ & $\Delta^2 h_0 a^3$ & $$ \\
$(96, 17, 48)$ & $\Delta^4 h_0$ & $$ \\
$(97, 17, 49)$ & $\Delta^4 h_1$ & $$ \\
$(99, 17, 50)$ & $\Delta^4 h_2$ & $$ \\
$(104, 19, 53)$ & $\Delta^4 c$ & $$ \\
$(104, 20, 52)$ & $\Delta^4 P$ & $$ \\
$(108, 22, 54)$ & $\Delta^4 h_0^3 a$ & $$ \\
$(110, 20, 56)$ & $\Delta^4 d$ & $$ \\
$(123, 22, 62)$ & $\Delta^4 a n$ & $$ \\
$(127, 24, 64)$ & $\tau^2 \Delta^4 d e$ & $\Delta^4 P^2 d$ \\
$(128, 23, 64)$ & $\tau^2 \Delta^4 a g$ & $$ \\
$(128, 23, 65)$ & $\Delta^4 (\Delta c + \tau a g)$ & $$ \\
$(132, 26, 66)$ & $\Delta^4 h_0 a^3$ & $$ \\
$(133, 24, 68)$ & $\tau^2 \Delta^4 e g$ & $\Delta^4 P d^2$ \\
$(140, 26, 70)$ & $\tau^2 \Delta^4 a^2 g$ & $\Delta^4 P^2 a n$ \\
$(144, 25, 72)$ & $\Delta^6 h_0$ & $\tau \Delta^4 P \Delta h_1 d$ \\
$(145, 27, 73)$ & $\tau^3 \Delta^4 a e g$ & $\Delta^4 P^2 \Delta c$ \\
$(146, 26, 74)$ & $\Delta^6 h_1^2$ & $\Delta^4 P a n d$ \\
$(147, 25, 74)$ & $\Delta^6 h_2$ & $$ \\
$(152, 27, 76)$ & $\tau \Delta^6 c$ & $$ \\
$(153, 28, 78)$ & $\Delta^6 h_1 c$ & $\tau^2 \Delta^4 P d^2 g$ \\
$(156, 30, 78)$ & $\Delta^6 h_0^3 a$ & $$ \\
$(168, 31, 84)$ & $\Delta^6 h_0 a^2$ & $$ \\
$(176, 31, 88)$ & $\tau \Delta^6 \Delta c$ & $$ \\
$(180, 34, 90)$ & $\Delta^6 h_0 a^3$ & $$ \\
$(192, 32, 96)$ & $\Delta^8$ & $$ \\
\end{longtable}

\begin{longtable}{ll}
\caption{$E_\infty$-page multiplicative generators
\label{tab:Einfty-gen}
} \\
\toprule
$(s,f,w)$ & $x$ \\
\midrule \endfirsthead
\caption[]{$E_\infty$-page multiplicative generators} \\
\toprule
$(s,f,w)$ & $x$ \\
\midrule \endhead
\bottomrule \endfoot
$(0, 1, 0)$ & $h_0$ \\
$(1, 1, 1)$ & $h_1$ \\
$(3, 1, 2)$ & $h_2$ \\
$(8, 3, 5)$ & $c$ \\
$(8, 4, 4)$ & $P$ \\
$(12, 6, 6)$ & $h_0^3 a$ \\
$(14, 4, 8)$ & $d$ \\
$(20, 4, 12)$ & $g$ \\
$(24, 7, 12)$ & $h_0 a^2$ \\
$(25, 5, 13)$ & $\Delta h_1$ \\
$(27, 6, 14)$ & $a n$ \\
$(32, 7, 16)$ & $\tau^2 a g$ \\
$(32, 7, 17)$ & $\Delta c + \tau a g$ \\
$(36, 10, 18)$ & $h_0 a^3$ \\
$(48, 10, 24)$ & $\Delta^2 h_0^2$ \\
$(51, 9, 26)$ & $\Delta^2 h_2$ \\
$(56, 11, 28)$ & $\tau \Delta^2 c$ \\
$(56, 13, 28)$ & $\Delta^2 P h_0 + \tau^3 \Delta h_1 d e$ \\
$(60, 14, 30)$ & $\Delta^2 h_0^3 a$ \\
$(72, 15, 36)$ & $\Delta^2 h_0 a^2$ \\
$(80, 15, 40)$ & $\tau \Delta^2 \Delta c$ \\
$(84, 18, 42)$ & $\Delta^2 h_0 a^3$ \\
$(96, 17, 48)$ & $\Delta^4 h_0$ \\
$(97, 17, 49)$ & $\Delta^4 h_1$ \\
$(99, 17, 50)$ & $\Delta^4 h_2$ \\
$(104, 19, 53)$ & $\Delta^4 c$ \\
$(104, 20, 52)$ & $\Delta^4 P$ \\
$(108, 22, 54)$ & $\Delta^4 h_0^3 a$ \\
$(110, 20, 56)$ & $\Delta^4 d$ \\
$(120, 23, 60)$ & $\Delta^4 h_0 a^2$ \\
$(123, 22, 62)$ & $\Delta^4 a n$ \\
$(128, 23, 65)$ & $\Delta^4 (\Delta c + \tau a g)$ \\
$(128, 23, 64)$ & $\tau^2 \Delta^4 a g$ \\
$(132, 26, 66)$ & $\Delta^4 h_0 a^3$ \\
$(144, 26, 72)$ & $\Delta^6 h_0^2$ \\
$(147, 25, 74)$ & $\Delta^6 h_2$ \\
$(152, 27, 76)$ & $\tau \Delta^6 c$ \\
$(152, 29, 76)$ & $\Delta^4 (\Delta^2 P h_0 + \tau^3 \Delta h_1 d e)$ \\
$(156, 30, 78)$ & $\Delta^6 h_0^3 a$ \\
$(168, 31, 84)$ & $\Delta^6 h_0 a^2$ \\
$(176, 31, 88)$ & $\tau \Delta^6 \Delta c$ \\
$(180, 34, 90)$ & $\Delta^6 h_0 a^3$ \\
$(192, 32, 96)$ & $\Delta^8$ \\
\end{longtable}

\begin{longtable}{lll}
\caption{Hidden $\tau$ extensions
\label{tab:hidden-tau}
} \\
\toprule
$(s,f,w)$ & source & target \\
\midrule \endfirsthead
\caption[]{Hidden $\tau$ extensions} \\
\toprule
$(s,f,w)$ & source & target \\
\midrule \endhead
\bottomrule \endfoot
$(22, 6, 14)$ & $h_1^2 g$ & $c d$ \\
$(22, 7, 13)$ & $c d$ & $P d$ \\
$(23, 7, 15)$ & $h_1^3 g$ & $h_1 c d$ \\
$(23, 8, 14)$ & $h_1 c d$ & $P h_1 d$ \\
$(28, 7, 17)$ & $c g$ & $d^2$ \\
$(29, 8, 18)$ & $h_1 c g$ & $h_1 d^2$ \\
$(35, 9, 19)$ & $\tau^2 h_1 d g$ & $P a n$ \\
$(40, 9, 23)$ & $\tau h_0 g^2$ & $\D h_1^2 d$ \\
$(40, 10, 22)$ & $\D h_1^2 d$ & $P (\D c + \tau a g)$ \\
$(41, 9, 23)$ & $\tau^2 h_1 g^2$ & $a n d$ \\
$(46, 10, 26)$ & $\D h_1^2 g$ & $(\D c + \tau a g) d$ \\
$(47, 11, 27)$ & $\D h_1^3 g$ & $\D h_1 c d$ \\
$(47, 12, 26)$ & $\D h_1 c d$ & $P \D h_1 d$ \\
$(53, 12, 30)$ & $\D h_1 c g$ & $\D h_1 d^2$ \\
$(71, 15, 39)$ & $\D^2 h_1^3 g$ & $\D^2 h_1 c d$ \\
$(118, 22, 62)$ & $\D^4 h_1^2 g$ & $\D^4 c d$ \\
$(118, 23, 61)$ & $\D^4 c d$ & $\D^4 P d$ \\
$(119, 23, 63)$ & $\D^4 h_1^3 g$ & $\D^4 h_1 c d$ \\
$(119, 24, 62)$ & $\D^4 h_1 c d$ & $\D^4 P h_1 d$ \\
$(124, 23, 65)$ & $\D^4 c g$ & $\D^4 d^2$ \\
$(125, 24, 66)$ & $\D^4 h_1 c g$ & $\D^4 h_1 d^2$ \\
$(131, 25, 67)$ & $\tau^2 \D^4 h_1 d g$ & $\D^4 P a n$ \\
$(136, 25, 71)$ & $\tau \D^4 h_0 g^2$ & $\D^5 h_1^2 d$ \\
$(136, 26, 70)$ & $\D^5 h_1^2 d$ & $\D^4 P (\D c + \tau a g)$ \\
$(137, 25, 71)$ & $\tau^2 \D^4 h_1 g^2$ & $\D^4 a n d$ \\
$(142, 26, 74)$ & $\D^5 h_1^2 g$ & $\D^4 (\D c + \tau a g) d$ \\
$(143, 27, 75)$ & $\D^5 h_1^3 g$ & $\D^5 h_1 c d$ \\
$(143, 28, 74)$ & $\D^5 h_1 c d$ & $\D^5 P h_1 d$ \\
$(149, 28, 78)$ & $\D^5 h_1 c g$ & $\D^5 h_1 d^2$ \\
$(150, 30, 78)$ & $\tau^8 \D^2 h_1^2 g^5$ & $\D^4 P (\D c + \tau a g) d$ \\
$(167, 31, 87)$ & $\D^6 h_1^3 g$ & $\D^6 h_1 c d$ \\
\end{longtable}

\section{Charts}

The following charts display
the $E_2$-page, $E_3$-page, $E_4$-page, and
$E_\infty$-page 
of the Adams spectral sequence for $\mmf$.  

The $E_2$-page is free as a module over $\F_2[\Delta^4]$, 
where $\Delta^4$ is a class, also known as $v_2^{16}$,
in the 96-stem.
The chart shows only
elements that are not multiples of $\Delta^4$.
Elements are labelled using the notation from \cite{Isaksen09}.

The $E_3$-page, $E_4$-page, and $E_\infty$-pages are free as modules
over $\F_2[\Delta^8]$, 
where $\Delta^8$ is a class, also known as $v_2^{32}$, in the 192-stem.  The charts show only
elements that are not multiples of $\Delta^8$.

The $E_\infty$-page is not free as a module over $\F_2[g]$ or
$\F_2[P]$, although it does contain some elements that are 
$g$-periodic and $P$-periodic.  One can easily identify
the $P$-periodic families as repeating patterns of slope $1/2$.

The $g$-periodic families are harder to identify.
Except for multiples of $P$, 
the element $\Delta^6 h_0 a^3$ in the $180$-stem is the last element
that is annihilated by $g$.  In other words, excluding multiples of $P$
and $\Delta^8$,
the $E_\infty$-page between the $201$-stem and the $220$-stem is the
same as the $E_\infty$-page between the $181$-stem and the $200$-stem.
No elements can be both $P$-periodic and $g$-periodic 
since $P^2 g$ is zero.

Each colored dot indicates a copy of $\M_2 = \F_2[\tau]$
or $\M_2/\tau^k$ for some $k \geq 1$.  
Table \ref{tab:dot-color} gives the interpretation of
each color.

Vertical lines, lines of slope $1$, and lines of slope $1/3$
indicate $h_0$ multiplications, $h_1$ multiplications, and
$h_2$ multiplications respectively.
Typically, the colors of these lines are determined by
the $\tau$ torsion of the target.
However,
magenta lines indicate that a product equals $\tau$ times a generator,
rather than a generator itself.
Also,
orange lines indicate that a product equals $\tau^k$ times a generator,
for some $k \geq 2$.  The exact value of $k$ can be deduced from the
weights of the source and target of the extension.

Red arrows of slope $1$ indicate infinite sequences of 
elements that are killed by $\tau$ and are connected by 
$h_1$-multiplications.

Light blue lines with negative slopes indicate differentials.
Lines of negative slope in other colors indicate that a differential
equals $\tau^k$ times a generator for some value of $k \geq 1$,
rather than a generator itself.  These colors are compatible
with the colors listed in Table \ref{tab:dot-color}.  
For example, a red line
indicates that a differential equals $\tau$ times a generator,
and a blue line indicates that a differential equals $\tau^2$ times
a generator.

In the $E_\infty$-page,
yellow lines indicate hidden $\tau$ extensions.

\begin{longtable}{ll}
\caption{Dot color interpretations
\label{tab:dot-color}
} \\
\toprule
\endfirsthead
\caption[]{Dot color interpretations} \\
\toprule
\endhead
\bottomrule \endfoot
$\M_2$ & gray \\
$\M_2/\tau$ & red \\
$\M_2/\tau^2$ & blue \\
$\M_2/\tau^3$ & green \\
$\M_2/\tau^4$ & pink \\
$\M_2/\tau^5$ & light green \\
$\M_2/\tau^6$ & dark blue \\
$\M_2/\tau^9$ & purple \\
$\M_2/\tau^{10}$ & lavender \\
$\M_2/\tau^{11}$ & brown \\
\end{longtable}

%\includepdf[lastpage=2]{charts/AdamsE5.pdf}

\newpage

\pdfpagewidth=108in

\pdfpageheight=57in

\includegraphics{AdamsE2.pdf}

\newpage

\includegraphics{AdamsE3.pdf}

\newpage

\includegraphics{AdamsE4.pdf}

\newpage

\includegraphics{AdamsE5.pdf}

\end{document}